\documentstyle[12pt]{article} 
\begin{document}
\newcommand{\singlespace}{
    \renewcommand{\baselinestretch}{1}
\large\normalsize}
\newcommand{\doublespace}{
   \renewcommand{\baselinestretch}{1.2}
   \large\normalsize}
\renewcommand{\theequation}{\thesection.\arabic{equation}}

\input amssym.def
\input amssym
\setcounter{equation}{0}
\def \ten#1{_{{}_{\scriptstyle#1}}}
\def \Z{\Bbb Z}
\def \C{\Bbb C}
\def \R{\Bbb R}
\def \Q{\Bbb Q}
\def \N{\Bbb N}
\def \l{\lambda}
\def \V{V^{\natural}}
\def \wt{{\rm wt}}
\def \tr{{\rm tr}}
\def \Res{{\rm Res}}
\def \End{{\rm End}}
\def \Aut{{\rm Aut}}
\def \mod{{\rm mod}}
\def \Hom{{\rm Hom}}
\def \im{{\rm im}}
\def \<{\langle} 
\def \>{\rangle} 
\def \w{\omega}
\def \c{{\tilde{c}}}
\def \o{\omega}
\def \t{\tau }
\def \ch{{\rm ch}}
\def \a{\alpha }
\def \b{\beta}
\def \e{\epsilon }
\def \la{\lambda }
\def \om{\omega }
\def \O{\Omega}
\def \qed{\mbox{ $\square$}}
\def \pf{\noindent {\bf Proof: \,}}
\def \voa{vertex operator algebra\ }
\def \voas{vertex operator algebras\ }
\def \p{\partial}
\def \1{{\bf 1}}
\def \ll{{\tilde{\lambda}}}
\def \H{{\bf H}}
\def \F{{\bf F}}
\def \h{{\frak h}}
\def \g{{\frak g}}
\def \rank{{\rm rank}}
\singlespace
\newtheorem{thmm}{Theorem}
\newtheorem{co}[thmm]{Corollary}
\newtheorem{th}{Theorem}[section]
\newtheorem{prop}[th]{Proposition}
\newtheorem{coro}[th]{Corollary}
\newtheorem{lem}[th]{Lemma}
\newtheorem{rem}[th]{Remark}
\newtheorem{de}[th]{Definition}
\newtheorem{con}[th]{Conjecture}
\newtheorem{ex}[th]{Example}

\begin{center}
{\Large {\bf  $W$-algebra $W(2,2)$ and the vertex operator algebra $L(\frac{1}{2},0)\otimes L(\frac{1}{2},0)$}} \\
\vspace{0.5cm}

Wei Zhang\  and \ Chongying Dong\footnote{Supported by NSF grants and a research grant from the
Committee on Research, UC Santa Cruz (dong@math.ucsc.edu).} \\
Department of mathematics\\
University of California\\
 Santa Cruz, CA 95064 

\end{center}
\hspace{1.5 cm}

\begin{abstract} 
In this paper the $W$-algebra $W(2,2)$ and its representation theory are 
studied. It is proved that a simple vertex operator algebra generated by 
two weight 2 vectors is either
a vertex operator algebra associated to a highest irreducible 
$W(2,2)$-module or a tensor product of two irreducible Virasoro vertex
operator algebras. Furthermore, any rational, $C_2$-cofinite and simple
vertex operator algebra whose weight 1 subspace is zero and weight 2 subspace
is 2-dimensional, and with central charge
$c=1$ is isomorphic to $L(\frac{1}{2},0)\otimes L(\frac{1}{2},0).$
\end{abstract}

\section{Introduction}

Motivated partially by the problem of classification of rational
vertex operator algebras with central charge $c=1$ and by the
Frenkel-Lepowsky-Meurman's uniqueness conjecture on the moonshine
vertex operator algebra $V^{\natural}$ \cite{FLM}, we give a
characterization of the vertex operator algebra $L(1/2,0)\otimes
L(1/2,0)$ in terms of the central charge and the dimensions of weights
1 and 2 subspaces in this paper.  Here the $L(1/2,0)$ is the
vertex operator algebra associated to the irreducible highest weight
module for the Virasoro algebra with central charge $1/2$ which is the
smallest central charge among the discrete unitary series for the
Virasoro algebra.

The classification of $c=1$ rational conformal field theories
at character level has been achieved in the physics literature
under the assumption that the sum of the square of the norm of
the irreducible characters is a modular function over the full modular
group \cite{K}. But the classification of rational vertex operator algebras
with $c=1$ remains an open and hard problem. If a vertex operator algebra   
$V=\sum_{n\geq 0}V_n$ with $\dim V_0=1$ is rational and $C_2$-cofinite,
then $V_1$ is a reductive Lie algebra  and its rank is less than
or equal to the effective central charge $\c$ \cite{DM2}. Also, the
vertex operator subalgebra generated by $V_1$ is a tensor product
of vertex operator algebras associated to integrable highest weight 
modules for affine Kac-Moody algebras and lattice vertex operator algebra
\cite{DM3}. In the case 
that $c=\c=1,$ we can classify the vertex operator algebras with
$\dim V_1\ne 0.$ Since $V_1$ is a reductive Lie algebra whose rank
is less than or eaual to 1, we immediately see that $v_1$ is either 
1-dimensional or 3-dimensional, as a result, $V$ is isomorphic to
a vertex operator algebra associated to a rank 1 lattice. So
one can assume that $V_1=0.$  
The vertex operator algebra $L(1/2,0)\otimes L(1/2,0)$ has this property.
So characterization of  $L(1/2,0)\otimes L(1/2,0)$ can be regarged as
a part of program of classification of rational vertex operator
algebras with $c=1.$

The vertex operator algebra  $L(1/2,0)\otimes L(1/2,0)$ plays an important
role in the study of the moonshine vertex operator algebra $V^{\natural}.$ 
The moonshine vertex operator algebra $V^{\natural}$ which is
fundamental in shaping the field of vertex operator algebra was
constructed as a bosonic orbifold theory based on the Leech lattice
\cite{FLM}.  The discovery of existence of $L(1/2,0)^{\otimes 48}$
inside the moonshine vertex operator algebra $V^{\natural}$
\cite{DMZ} opens a different way to study $V^{\natural}.$ This
leads to  the theory of code and framed vertex operator algebras
\cite{M2},\cite{DGH}. This discovery is also essential in 
a proof that $V^{\natural}$ is holomorphic
\cite{D}, a new construction of $V^{\natural}$ \cite{M3},
proofs of weak versions of the Frenkel-Lepowsky-Meurman's
uniqueness conjecture on $V^{\natural}$ \cite{DGL}, \cite{LY} and
a study of $V^{\natural}$ in terms of conformal nets \cite{KL}. There is
no doubt that a characterization of $L(1/2,0)\otimes L(1/2,0)$ will be
very helpful in the study of structure of $V^{\natural}$ and the 
 Frenkel-Lepowsky-Meurman's
uniqueness conjecture.

The $W(2,2)$ and its highest weight modules enter the picture naturally
during our discussion on $L(1/2,0)\otimes L(1/2,0).$ The $W$-algebra
$W(2,2)$ is an extension of the Virasoro algebra and also has a very good
highest weight module theory (see Section 2). Its highest weight modules
produce a new class of vertex  operator algebras. Contrast to the 
Virasoro algebra case, this class of vertex operator algebras
are always irrational.
From this point of view, this class of vertex operator algebras are not
interesting.

The $W(2,2)$ and associated vertex operator algebras are also closely
related to the classification of simple vertex operator algebra with 2
generators. It is well known that each homogeneous subspace $V_n$ of a
vertex operator algebra $V=\sum_{n\in\Z}V_n$ is some kind of algebra
under the product $u\cdot v=u_{n-1}v$ for $u,v\in V_n$ where $u_m$ is
the component operator of $Y(u,z)=\sum_{m\in \Z}u_mz^{-m-1}.$ If a
vertex operator algebra $V=\sum_{n\geq 0}V_n$ with $\dim V_0=1$ is
rational and $C_2$-cofinite, then $V_1$ and the vertex operator
subalgebra generated by $V_1$ are well understood \cite{DM2}. So it is
natural to turn our attention to $V_2.$ This is still a very hard
problem even with $V_1=0.$ A simple vertex operator algebra $V$
satisfying $V_1=0$ is called the {\em moonshine type}. The $V_2$ in
this case is a commutative nonassociative algebra. The simple vertex
operator algebras of the moonshine type with $\dim V_2=2$ and generated by
$V_2$ are also
classified in this paper.  There are two families of such
algebras. One of this family consists the tensor product of two vertex
operator algebras associated to the irreducible highest weight
modules for the Virasoro algebra and the other family consists the
vertex operator algebras associated the highest weight modules for the
the $W$-algebra $W(2,2).$

The paper is organized as follows. We define and study the $W$-algebra
$W(2,2)$ in Section 2. In particular we use the bilinear form on Verma
modules $V(c,h_1,h_2)$ to determine the irreducible quotient modules
$L(c,h_1,h_2)$ for $W(2,2)$ for most $c$ and $h_i.$ In Section 3 we classify
the simple vertex operator algebras of the moonshine type generated by 2 weight
2 vectors. The section 4 is devoted to the characterization of 
rational vertex operator algebra $L(1/2,0)\otimes L(1/2,0).$ The main idea
is to use the modular invariance of the graded characters of the irreducible
modules \cite{Z} to control the growth of the graded dimensions of the vertex
operator algebra.  

\section{$W$-algebra $W(2,2)$}
\setcounter{equation}{0}

The $W$-algebra $W(2,2)$ considered in this paper is an infinite dimensional
Lie algebra with generators $L_m, W_m,C$ for $m\in \Z$ and Lie bracket
$$[L_m,L_n]=(m-n)L_{m+n}+\frac{m^3-m}{12}\delta_{m+n,0}C,$$
$$[L_m,W_n]=(m-n)W_{m+n}+\frac{m^3-m}{12}\delta_{m+n,0}C,$$
$$[W_m,W_n]=0$$
for $m,n\in \Z$ where $C$ is a central element. In this section we study
the highest weight modules for this algebra and the corresponding 
vertex operator algebras. 

Let $c,h_1,h_2\in\C$ and we denote by $V(c,h_1,h_2)$ the highest
weight module for $W(2,2)$ with central charge $c$ and highest weight
$(h_1,h_2).$ Then  $V(c,h_1,h_2)=U(W(2,2))/I_{c,h_1,h_2}$ where
$I_{c,h_1,h_2}$ is the left ideal of the universal enveloping algebra
$U(W(2,2))$ generated by $L_m,W_m,$ $C-c,$ $L_0-h_1$ and $W_0-h_2$ for
positive $m.$ The $V(c,h_1,h_2)$ can also be realized as induced module
as in the case of Virasoro algebra. It is standard that $V(c,h_1,h_2)$
has a unique maximal submodule $J(c,h_1,h_2)$ so that $L(c,h_1,h_2)=
V(c,h_1,h_2)/J(c,h_1,h_2)$ is an irreducible highest weight module. 
As in the case of Virasoro algebra, there is a unique invariant symmetric
bilinear form $(,)$ on $V(c,h_1,h_2)$  such that
$$(L_mu,v)=(u,L_{-m}v),\ \  (W_mu,v)=(u,W_{-m}v),$$
$$(1,1)=1$$
where $1=1+I_{c,h_1,h_2}.$ Moreover, the radical of this
bilinear form is exactly the maximal submodule $J(c,h_1,h_2).$

Let $X$ be a proper submodule of $V(c,h_1,h_2).$ Then $X$ is a submodule
of $J(c,h_1,h_2)$ and the bilinear form  $(,)$ on  $V(c,h_1,h_2)$ induces
a symmetric invariant bilinear from $(,)$ on the quotient module
 $V(c,h_1,h_2)/X.$ 

As in the classical case we need to answer the basic question: What is
 $J(c,h_1,h_2)?$  We first consider the case $(c,h_1,h_2)=(c,0,0).$ Clearly,
$L(0,0,0)=\C.$  So we now assume that $c\ne 0.$ Note that $U(W(2,2))L_{-1}1
+U(W(2,2))W_{-1}1$ is a proper submodule of $V(c,0,0).$  

\begin{th}\label{t2.1} If $c\ne 0$ then $J(c,0,0)= U(W(2,2))L_{-1}1
+U(W(2,2))W_{-1}1$ and $L(c,0,0)$ has basis
$$S=\{W_{-m_1}\cdots W_{-m_s}L_{-n_1}\cdots L_{-n_t}\1|m_1\geq \cdots \geq m_s>1,
n_1\geq \cdots \geq  n_t>1\}$$
where $\1$ is the canonical highest weight vector of $L(c,0,0).$ 
\end{th}

\pf. Set $\bar{V}(c,0,0)=V(c,0,0)/(U(W(2,2))L_{-1}1+U(W(2,2))W_{-1}1.$ Then
$S$ forms a basis of $\bar{V}(c,0,0)$ by PBW theorem.  
Let $S_n$ be the subset of $S$ consisting of weight $n$ (with respect to 
the operator $L_0$) vectors:
$$W_{-m_1}\cdots W_{-m_s}L_{-n_1}\cdots L_{-n_t}\1$$
with $m_1\geq \cdots \geq m_s>1,
n_1\geq \cdots \geq n_t>1$ and $\sum_im_i+\sum_jn_j=n.$
So it is enough to show that $S_n$ is linearly independent in $L(c,0,0).$  
The idea is to prove that the determinant of $((u,v))_{u,v\in S_n}$ is nonzero.

To see how the determinant argument works we first consider two subsets
$S_{0,n}$ and $S_{n,0}$ of $S_n$ where $S_{0,n}$ consists of vectors 
$L_{-n_1}\cdots L_{-n_t}\1$ for $n_1\geq \cdots \geq \cdots n_t>1$ 
and $S_{n,0}$ consists of vectors 
$W_{-m_1}\cdots W_{-m_s}\1$ for positive integers $m_1\geq \cdots \geq m_s>1.$

Let
$$P_n=\{(m_1,m_2\cdots m_s)|s\geq 1,m_1\geq m_2\geq\cdots \geq m_s>1,\sum_{i}m_i=n\}$$
which is a set of partitions of $n$ without $1.$ 
We define a total order on $P_n$ so that if $(m_1,...,m_s)>(p_1,...,p_t)$
if and only if there exists $1\leq k\leq s$ such that 
$m_i=p_i$ for $i<k$ and $m_k>n_k.$ Note that each vector in $S_{0,n}$ or $S_{n,0}$ 
is associated to a partition in $P_n.$ 
We can label the corresponding vectors in $A$ and $B$ 
by using this order. Let $p_n=|P_n|$ be the cardinality of $P_n.$ Then we can write the
vectors in $S_{0,n}$ by $\{u_1^{0,n},...,u_{p_n}^{0,n}\}$ so that if $i<j$ then 
the partition corresponding to $u_i^{0,n}$ is less than the partition corresponding
to $u_j^{0,n}.$ We denote the elements in $B$ by $\{v_1^{n,0},...,v_{p_n}^{n,0}\}$ with
the same order.
 
Observe that if $(m_1,...,m_s)\in P_n$
and $m\geq m_1$ then 
$$L_{m}W_{-m_1}\cdots W_{-m_s}\1= \frac{m^3-m}{12}c\frac{\partial}{\partial W_{-m}}W_{-m_1}\cdots W_{-m_s}\1.$$
This implies immediately that if $i>j$ $(u_i^n,v_j^n)=0.$ Clearly, 
$(u_i^n,v_i^n)\ne 0.$ This shows that both $S_{0,n}$ and $S_{n,0}$ are linearly
independent in $L(c,0,0).$ 

Let
$$u=W_{-m_1}\cdots W_{-m_s}L_{-n_1}\cdots L_{-n_t}\1,$$
$$v=W_{-p_1}\cdots W_{-p_a}L_{-q_1}\cdots L_{-q_b}\1\in S_n.$$ 
Clearly, $(u,v)=0$ if $\sum m_i>\sum q_j$ or $\sum p_i>\sum n_j.$ 
If  $\sum m_i=\sum q_j$ and $\sum n_i=\sum p_j$ we have
$$(u,v)=(W_{-m_1}\cdots W_{-m_s}\1,L_{-q_1}\cdots L_{-q_b}\1)\cdot$$
$$(W_{-p_1}\cdots W_{-p_a}\1, L_{-n_1}\cdots L_{-n_t}\1).$$ 
which follows from the computation
\begin{eqnarray*}
& &\ \ \ (W_{-m_1}\cdots W_{-m_s}\1,L_{-q_1}\cdots L_{-q_b}\1)\1\\
& &=L_{q_b}\cdots L_{q_1}W_{-m_1}\cdots W_{-m_s}\1\\
& &=W_{m_s}\cdots W_{m_1}L_{-q_1}\cdots L_{-q_b}\1,
\end{eqnarray*}
where we have used the fact that $W_m$ commute with each other.

We now fix nonnegative integers $d$ and $e$ such that $d+e=n.$ Let
$(m_1,...,m_s)\in P_d$ and $(n_1,...,n_t)\in P_e.$ Assume that 
$u_i^d=L_{-m_1}\cdots L_{-m_s}\1$ and $u_j^e=L_{-n_1}\cdots L_{-n_t}\1$
(see the argument above for the definition of $u_i^n$). Set
$$u_{i,j}^{d,e}=W_{-m_1}\cdots W_{-m_s}L_{-n_1}\cdots L_{-n_t}\1$$
for $i=1,...,p_d$ and $j=1,...,p_e$ and let $S_{d,e}$ the set consisting
of these vectors.  Then $S_{n-1,1}=S_{1,n-1}=\emptyset$  and 
$$S_n=\cup_{d\ne 1,n-1}S_{d,n-d}.$$
From the discussion above, for any $u\in S_{d,e}$ and $v\in S_{d_1,e_1}$
with $d,d_1\ne 1$ we have $(u,v)=0$ if $d>e_1$ or $d_1>e.$ 

Again from the discussion on $S_{0,n}$ and $S_{n,0}$ we can relabel the vectors in $S_{d,n-d}$ 
in two different ways such as $\{u^{d,n-d}_{i}|i=1,...,p_dp_e\}$
and  $\{v^{d,n-d}_{i}|i=1,...,p_dp_e\}$ so that 
$(u^{d,n-d}_i,v^{n-d,d}_j)=0$ if $i>j$ and $(u^{d,n-d}_i,v^{n-d,d}_i)$ is 
nonzero for all $i.$

Let 
$$A=\left(\begin{array}{ccccccc}
A_{0,0} & A_{0,2}, &A_{0,3}& \cdots & A_{0,n-3} & A_{0,n-2} & A_{0,n}\\ 
A_{2,0} & A_{2,2}, &A_{2,3}& \cdots & A_{2,n-3} & A_{2,n-2} & A_{2,n}\\
A_{3,0} & A_{2,2}, &A_{3,3}& \cdots & A_{3,n-3} & A_{3,n-2} & A_{3,n}\\ 
\vdots  & \vdots  & \vdots  & \vdots  & \vdots  & \vdots  & \vdots  \\
A_{n-3,0} & A_{n-3,2}, &A_{n-3,3}& \cdots & A_{n-3,n-3} & A_{n-3,n-2} & A_{n-3,n}\\
A_{n-2,0} & A_{n-2,2}, &A_{n-2,3}& \cdots & A_{n-2,n-3} & A_{n-2,n-2} & A_{n-2,n}\\
A_{n,0} & A_{n,2}, &A_{n,3}& \cdots & A_{n,n-3} & A_{n,n-2} & A_{n,n}
\end{array}\right)$$
where  $A_{d,e}$ is a submatrix defined by 
$$A_{d,e}=((u^{d,n-d}_i,v^{n-e,e}_j))$$
with obvious ranges for $i,j.$ Clearly, $A_{d,e}=0$ if $d>e.$ So $A$ is 
an upper triangular matrix with nonzero diagonals.  Thus $A$ 
is nondegenerate. As a result 
$S_n$ is linearly independent. This finishes the proof.
\qed

\begin{rem}{\rm  Although $W(2,2)$ is an extension of the Virasoro algebra,
the representation theory for $W(2,2)$ is different from that for
the Virasoro algebra in a fundamental way. For $W(2,2),$ the structure of
$L(c,0,0)$ for $c\ne 0$ is uniform and simple. But for the Virasoro
algebra, the situation is totally different. Let $L(c,h)$ be the 
irreducible highest weight module for the Virasoro algebra with central charge
$c$ and highest weight $h.$ In the case $c\ne 1-6(p-q)^2/pq$ where $p,q$ are
two coprime positive integers $1<p<q,$ then $L(c,0)=\bar{V}(c,0)$ where
$\bar{V}(c,0)=V(c,0)/U(Vir)L_{-1}v$ and $v$ is a nonzero highest weight vector
of the Verma module $V(c,0)$ (see \cite{FF}). The structure of $L(c_{s,t},0)$ 
is much more complicated. On the other hand, from the point of view of
vertex operator algebra, $L(c_{s,t},0)$ is a rational vertex operator 
algebra for all $c_{s,t}$ but $L(c,0)$ is not if $c\ne c_{s,t}$ 
(see \cite{FZ} and \cite{W}).}
\end{rem}

Next we discuss the vertex operator algebras associated to the highest weight 
modules for $W(2,2).$ Let $\1$ be  the canonical highest weight vector
of $V(c,0,0).$ From the axiom of vertex operator algebra we must modulo out
the submodule generated by $L_{-1}\1.$ From the commutator relation between
$L_m$ and $W_n$ we also know that $W_n$ cannot be the component operators
of a vertex operator associated to a weight one vector. This forces 
$W_{-1}\1=0$ if there is a vertex operator algebra structure. So by Theorem
\ref{t2.1}, $L(c,0,0)$ is the only quotient of $V(c,0,0)$ which may have a structure of vertex operator algebra.

A $W(2,2)$-module $M$ is restricted if for any $w\in M,$
$L_{m}w=W_{m}w=0$ if $m$ is sufficiently large. Recall the weak module, 
admissible module and
ordinary module from \cite{DLM1}.

\begin{th}\label{t2.3} Assume that $c\ne 0.$ Then

(1) There is a unique vertex operator
algebra structure on $L(c,0,0)$ with vacumme $\1$ and Virasoro element
$\omega=L_{-2}\1.$ Moreover, $L(c,0,0)$ is generated by $\omega$ and
$x=W_{-2}\1$ such that $Y(\omega,z)=\sum_{n\in\Z}L_nz^{-n-2}$
and $Y(x,z)=\sum_{n\in\Z}W_nz^{-n-2}.$

(2) If $M$ is a restricted $W(2,2)$-module with central charge $c,$ then $M$ is
a weak $L(c,0,0)$-module such that $Y_M(\omega,z)=\sum_{n\in\Z}L_nz^{-n-2}$
and $Y_M(x,z)=\sum_{n\in\Z}W_nz^{-n-2}.$ In particular, any quotient 
module of $V(c,h_1,h_2)$ is an ordinary module for $L(c,0,0).$

(3) Any irreducible admissible $L(c,0,0)$-module is ordinary.

(4) $\{L(c,h_1,h_2)|h_i\in \C\}$ gives a complete list of irreducible
$L(c,0,0)$-modules up to isomorphism.
\end{th}  

\pf (1) and (2) are fairly standard following from the local system theory
(see \cite{L2}, \cite{LL}). (3) and (4) follow from that fact that 
any irreducible admissible module for $L(c,0,0)$ is an irreducible highest 
weight module for $W(2,2).$ \qed
\bigskip

We now turn our attention to the Verma module $V(c,h_1,h_2)$ in general.
As in general highest weight module theory, we want to know 
when  $V(c,h_1,h_2)=L(c,h_1,h_2)$
is irreducible.
\begin{th} The Verma module $V(c,h_1,h_2)$ is irreducible if and only if
$\frac{m^2-1}{12}c+2h_2\ne 0$ for any nonzero integer $m.$ 
\end{th}

\pf As in the proof of Theorem \ref{t2.1} we use the determinant of the 
invariant bilinear form to prove the result. Note that 
$$V(c,h_1,h_2)=\oplus_{n\geq 0}V(c,h_1,h_2)_{h_1+n}$$
where $V(c,h_1,h_2)_n$ has a basis consisting of vectors
$$W_{-m_1}\cdots W_{-m_s}L_{-n_1}\cdots L_{-n_t}\1$$
where $m_1\geq \cdots \geq m_s>0,
n_1\geq \cdots \geq  n_t>0, \sum{m_i}+\sum{n_j}=n.$
As in the proof of Theorem \ref{t2.1} we still denote this set by $S_n.$
Define a matrix $A_n=((u,v))_{u,v\in S_n}.$  Then $V(c,h_1,h_2)=L(c,h_1,h_2)$
if and only if $\det A_n\ne 0$ for all $n>0.$

Note that if $m\geq m_1\geq \cdots\geq m_s>0$  
$$L_{m}W_{-m_1}\cdots W_{-m_s}\1= (\frac{m^3-m}{12}c+2mh_2)\frac{\partial}{\partial W_{-m}}W_{-m_1}\cdots W_{-m_s}\1.$$ 
Following the proof of Theorem \ref{t2.1} we see immediately that
if $\frac{m^2-1}{12}c+2h_2\ne 0$ for all $0\ne m\in\Z$ then $\det A_n\ne 0$
and if $\frac{m^2-1}{12}c+2h_2= 0$ for some $0<m,$ then 
$\det A_m=0.$ The proof is complete. \qed

\bigskip

It is definitely interesting to determine the $J(c,h_1,h_2)$ if 
$\frac{m^2-1}{12}c+2h_2=0$ for some nonzero integer $m.$ But this will
be a problem which has nothing to do with  the characterization of 
$L(1/2,0)\otimes L(1/2,0)$ in this paper. 
 We will not go in this direction further.

\section{Vertex operator algebras of the moonshine type}
\setcounter{equation}{0}

In this section we classify the simple vertex operator algebras $V$ of
the moonshine type such that $V$ is generated by $V_2$ and $V_2$ is
2-dimensional.

Note that $V_0=\C\1$ is 1-dimensional for the moonshine type vertex
operator algebra $V$ and $V_n=0$ if $n<0$ \cite{DGL}. Since $V_1=0$
and $V_0$ is 1-dimensional, there is a unique symmetric, nondegenerate
invariant bilinear from $(,)$ on $V$ such that $(\1,\1)=1$ (see
\cite{L1}). Then for any $u,v,w\in V$
$$(Y(u,z)v,w)=(v,Y(e^{L(1)z}(-z^{-2})^{L(0)}u,z^{-1})w)$$
and 
$$(u,v)\1=\Res_zz^{-1} Y(e^{L(1)z}(-z^{-2})^{L(0)}u,z^{-1})v.$$
In particular, the restriction of the form to each homegeneous 
subspace $V_n$ is nondegenerate and 
$$(u_{n+1}v,w)=(v,u_{-n+1}w)$$ for
all $u,v\in V_2$ and $w\in V.$
 
The $V_2$  is a commutative and associative 
algebra with the product $ab=a_1b$ for $a,b\in V_2$ and identity $\frac{\omega}{2}$ (cf. \cite{FLM})
The $V_2$ is called  Griess algebra of $V.$ 
Note that for $a,b\in V_2$ we have $(a,b)=a_3b.$ Moreover, the form
on $V_2$ is associative. That is, $(ab,c)=(a,bc)$ for $a,b,c\in V_2.$ 

\begin{th}\label{t3.1} Let $V$ be a simple vertex operator algebra of
the moonshine type with central charge $c\ne 0$ 
such that $V$ is generated by $V_2$ and $V_2$ is
2-dimensional. Then $V$ is isomorphic to $L(c_1,0)\otimes L(c_2,0)$ for
some nonzero complex number $c_1,c_2$ such that $c_1+c_2=c$ if
$V_2$ is semisimple, and isomorphic to $L(c,0,0)$ if $V_2$ is not semisimple.
\end{th}

\pf  If $V_2$ is a 2-dimensional semisimple commutative associative
algebra with the identity $\omega/2.$ Then $\omega=\omega^1+ \omega^2$
so that $\omega^1/2$ and $\omega^2/2$ are the primitive idempotents. It follows
from \cite{M1} that $\omega^1$ and $\omega^2$ are Virasoro vectors. 
Let 
$$Y(\omega^i,z)=\sum_{n\in\Z}L^i(n)z^{-n-2}$$
for $i=1,2.$
Then 
$$[L^i(m),L^i(n)]=(m-n)L^i(m+n)+\frac{m^3-m}{12}\delta_{m+n,0}c_i$$
for all $m,n\in \Z$ where $c_i\in \C$ is the central charge of 
$\omega^i.$ Since $\frac{\omega^i}{2}\frac{\omega^j}{2}=\delta_{i,j}$
we see that $(\omega^1)_3\omega^2=(\omega^1,\omega^2)=0$ by using the 
invariant property of the bilinear form. This implies that
$$[L^1(m),L^2(n)]=0$$
for all $m,n\in \Z$ and $c_1+c_2=c.$ 
Then $V=\<\omega^1\>\otimes \<\omega^2\>$ 
where $\<\omega^i\>$ is the vertex operator subalgebra of $V$ generated
by $\omega^i$ (with a different Virasoro vector). Note that 
$\<\omega^i\>$ is a quotient of $\bar{V}(c_i,0).$ Since $V$ is simple
we immediately have that $\<\omega^i\>$ is isomorphic to $L(c_i,0).$
As a result, $V$ is isomorphic to $L(c_1,0)\otimes L(c_2,0)$ in this case.

It remains to deal with the case that $V_2$ is not semisimple. Then
the Jacobson radical $J$ of $V_2$  is 1-dimensional. Assume that $J=\C x.$ 
Then $x^2=0$ and $(x,x)=(\omega/2, x^2)=0.$ Using the skew symmetry 
$Y(x,z)x=e^{L(-1)z}Y(x,-z)x$ we see that 
$$x_0x=-x_0x+L(-1)x_1x=-x_0x+L(-1)x^2=-x_0x.$$ 
This implies $x_0x=0.$ As a consequence, we see the component operators
$x_n$ of $Y(x,z)$ commute with each other. That is,
$[x_m,x_n]=0$ for all $m,n\in \Z.$ 

 Note that $(\omega,\omega)=L(2)\omega=c.$ Since the 
form $(,)$ on $V_2$ is nondegenerate, we may choose $x$ so that $(\omega,x)=c.$ Set
$W(m)=x_{m+1}$ for $m\in\Z.$  
Then we have the following commutator formula
\begin{eqnarray*}
[L(m),W(n)]=(m-n)W(m+n)+\frac{m^3-m}{12}\delta_{m,-n}c.
\end{eqnarray*}
This exactly says that the operators $L(m), W(m), c$ generate a copy
of $W(2,2)$ and $V$ is an irreducible highest weight module 
for $W(2,2).$ So in this case, $V$ is isomorphic to $L(c,0,0),$
as desired. \qed
\begin{rem}{\rm Theorem \ref{t3.1} is the main reason we introduce and
study the Lie algebra $W(2,2)$ and its highest weight modules. The vertex
operator algebra $L(c,0,0)$ will be used in the next section
when we characterize the rational vertex operator algebra $L(1/2,0)\otimes L(1/2,0).$}
\end{rem}

\section{Characterization of $L(1/2,0)\otimes L(1/2,0)$}
\setcounter{equation}{0}

In this section we give a characterization for the vertex operator algebra
$L(1/2,0)\otimes L(1/2,0).$ 

We first recall some basic facts about a rational vertex operator algebra
following from \cite{DLM1}. A vertex operator algebra $V$ is called rational
if any admissible module is completely reducible. It is proved in \cite{DLM1}
(also see \cite{Z}) that if $V$ is rational then there are only finitely
many irreducible admissible modules $M^1,...,M^k$ up to isomorphism
such that 
$$M^i=\oplus_{n\geq 0}M^i_{\l_i+n}$$
where $\l_i\in\Q,$  $M^i_{\l_i}\ne 0$ and each $M^i_{\l_i+n}$
is finite dimensional (see \cite{AM} and \cite{DLM2}). Let
$\l_{min}$ be the minimum of $\l_i$'s. The effective central charge
$\c$ is defined as $c-24\l_{min}.$ A vertex operator algebra is called
$C_2$-cofinite if $C_2(V)$ has finite codimension where 
$C_2(V)=\<u_{-2}v|u,v\in V\>.$ 

For each $M^i$ we define the $q$-character of $M^i$ by
$$\ch_q M^i=q^{-c/24}\sum_{n\geq 0}(\dim M^i_{\l_i+n})q^{n+\l_i}.$$
Then $\ch_q M^i$ converges to a holomorphic function on the upper
half plane if $V$ is $C_2$-cofinite \cite{Z}. Using the modular invariance
result from \cite{Z} and results on vector valued modular forms 
from \cite{KM} we have (see \cite{DM2})
\begin{lem}\label{growth} Let $V$ be rational and $C_2$-cofinite. For each $i,$
the coefficients of $\eta(q)^{\c}\ch_qM^i$ satisfy the polynomial
growth condition where
$$\eta(q)=q^{1/24}\prod_{n\geq 1}(1-q^n).$$
\end{lem}

We also need  some basic facts about the highest weight modules for the 
Virasoro algebra (see \cite{FF}, \cite{FQS},\cite{GKO}, \cite{FZ}, \cite{W}). 

\begin{prop}\label{vir} Let $c$ be a complex number. 

(i). $\bar{V}(c,0)$ is a vertex operator algebra and $L(c,0)$ is a simple
vertex operator algebra. 

(ii) The following are equivalent: (a) $\bar{V}(c, 0) = L(c, 0),$ 
(b) $c\ne c_{s,t}=1-6(s-t)^2/st$ for all coprime 
positive integers $s,t$ with $1<s<t,$ 
(c)  $L(c,0)$ is not rational. 
In this case, the $q$-graded character of $L(c,0)$ is equal 
to $\frac{q^{-c/24}}{\prod_{n>1}(1-q^n)}$ and the coefficients grow faster than
any polynomials. 

(iii) The following are equivalent: (a) $\bar{V}(c, 0) \ne L(c, 0),$ 
(b) $c=c_{s,t}$ for some $s,t,$ (c) $L(c,0)$ is rational.  
\end{prop}

From now on we assume that $V$ is a rational and $C_2$-cofinite vertex operator
algebra of the moonshine type such that $c=\c=1$ and $\dim V_2=2.$ We have already mentioned in Section 3 that $V_2$ is a commutative 
algebra with identity $\frac{\omega}{2}.$ The assumption that $\dim V_2=2$
makes $V_2$ an associative algebra. 

\begin{lem}\label{l4.3}  The $V_2$ is a semisimple associative algebra. That is, $V_2$ is
a direct sum of two ideals isomorphic to $\C.$
\end{lem}

\pf Suppose that $V_2$ is not semisimple. Recall from the proof of Theorem 
\ref{t3.1} that the Jacobson radical $J=\C x$ is one-dimensional. 
We assume that $(\omega,x)=1.$ Then the component operator $W(n)$ of $Y(x,z)=\sum_{n\in\Z}W(n)z^{-n-2}$ and the component operator of the $Y(\omega,z)$ generate a copy of the
$W$-algebra $W(2,2)$ with central charge $1.$ 

Let $U$ be the vertex operator subalgebra of $V$ generated by $V_2.$  Then
$U$ is isomorphic to $L(1,0,0)$ as a module for $W(2,2).$ 
From Theorem \ref{t2.1}, 
$$\ch_qU=\frac{q^{-1/24}}{\prod_{n>1}(1-q^n)^2}.$$
Note that $\ch_q U\leq \ch_qV,$ that is, the   coefficients of $\ch_qU$
are less than or equal to
the corresponding coefficients of $\ch_q V$. Then  $\ch_q U\leq \ch_qV$
and $\eta(q)\ch_qU\leq  \eta(q)\ch_qV$ as functions for $q\in (0,1).$
By Lemma  \ref{growth}, the coefficients of $\eta(q)\ch_qV$ satisfy the polynomial 
growth condition. On the other hand, the coefficients of
 $\eta(q)\ch_qU=\frac{1-q}{\prod_{n>1}(1-q^n)}$ grow faster than 
any polynomial in $n.$  This is a contradiction. 
\qed
\bigskip

So from now on we assume that $V_2$ is semisimple.
Again from the proof of Theorem \ref{t3.1}, we can write 
$\omega=\omega^1\oplus \omega^2$
so that $\omega^1/2$ and $\omega^2/2$ are the primitive idempotents.
The $\omega^1$ and $\omega^2$ are Virasoro vectors with central charges
$c_1$ and $c_2$ such that $c_1+c_2=1.$ Let $L^i(n)$ be as in Section 3. Then we have
two commutative Virasoro algebras: 
$$[L^i(m),L^j(n)]=\delta_{i,j}((m-n)L^i(m+n)+\frac{m^3-m}{12}\delta_{m+n,0}c_i)$$
for all $m,n\in \Z$ for $i,j=1,2.$ As before 
 we denote by $U$ the vertex operator subalgebra
of $V$ generated by $U.$ Then $U=\<\omega^1\>\otimes \<\omega^2\>$ 
where $\<\omega^i\>$ is the vertex operator subalgebra of $V$ generated
by $\omega^i$ (with a different Virasoro vector). Then
$\<\omega^i\>$ is a quotient of $\bar{V}(c_i,0).$ 

\begin{lem}\label{l4.2} If $c\ne 0,$ then the coefficients of $\ch_qL(c,0)$ does not 
satisfy the polynomial growth condition.
\end{lem}

\pf If $c\ne c_{s,t}$ for any coprime integers $1<s<t$ 
 then $\ch_qL(c,0)=\frac{q^{-c/24}}{\prod_{n>1}(1-q^n)}$
by Proposition \ref{vir} and the result is clear. We now assume that 
$c=c_{s,t}$ for some $s,t.$ Suppose that the coefficients of 
$$\ch_qL(c_{s,t},0)=q^{-c/24}\sum_{n\geq 0}a_nq^n$$ 
satisfy the polynomial growth condition. Then there exists a positive integer
$A$ and $\alpha$ such that $a_n\leq An^{\alpha}$ for all $n\geq 0.$

Let $m$ be a positive integer such that $m\geq \alpha.$ Then
$$\frac{1}{(1-q)^{m+1}}=\sum_{n\geq 0}{-m-1\choose n}(-1)^nq^n$$
where 
$${-m-1\choose n}=\frac{(-m-1)(-m-2)\cdots (-m-n)}{n!}={m+n\choose m}(-1)^n.$$
Thus 
$$\frac{1}{(1-q)^{m+1}}=\sum_{n\geq 0}{m+n\choose m}q^n.$$
Since ${m+n\choose m}$ is greater than $\frac{n^m}{m!}$ we see that
$$q^{c/24}\ch_qL(c,0)\leq m!A \frac{1}{(1-q)^{m+1}}$$
as formal power series.

We next prove that there exists a positive integer $k$ such that 
$kc_{s,t}\ne c_{s_1,t_1}$ for any coprime integers  $1<s_1<t_1.$ To see this we
need to examine the equation
$$1-\frac{6(s_1-t_1)^2}{st}=k(1-\frac{6(s-t)^2}{st})$$
which is equivalent to 
$$st(13s_1t_1-6s_1^2-6t_1^2)=s_1t_1k(13st-6s^2-6t^2).$$
Then both $s_1$ and $t_1$ are factors of $6st.$ So there are only finitely many
$s_1,t_1$ satisfy this equation. This implies that such $k$ exists.

Consider vertex operator algebra $L(c,0)^{\otimes k}$ which  contains
a vertex operator subalgebra $\bar{V}(ck,0)=L(ck,0)$ as $ck\ne c_{s_1,t_1}$ 
for any $s_1,t_1.$ So
$$q^{ck/24}\ch_q L(ck,0)\leq q^{ck/24}(\ch_q L(c,0)^{\otimes k})
=q^{ck/24}(\ch_q L(c,0))^k\leq (m!A)^k\frac{1}{(1-q)^{mk}}$$
and the coefficients of $q^{ck/24}\ch_q L(ck,0)$ satisfy the polynomial
growth condition.  

On the other hand we know from Proposition \ref{vir} that
$$q^{ck/24}\ch_q L(ck,0)=\frac{1}{\prod_{n>1}(1-q^n)}$$
whose coefficients satisfy the exponential growth condition. This is 
a contradiction. The proof is complete. \qed

\begin{lem}\label{l4.4} Let $\omega^i$ and $c_i$ be as before. Then 
$c_i=c_{s_i,t_i}$ for some coprime integers $1<s_i<t_i$ and
$\<\omega^i\>$ is isomorphic to $L(c_{s_i,t_i},0)$ for $i=1,2.$
\end{lem}

\pf First we note that as formal power series, $\ch_qU\leq \ch_qV.$
Let $U^i=\<\omega^i\>.$ Then $U=U^1\otimes U^2$ 
and $\ch_qU^1 ch_qU^2\leq \ch_qV.$ 
Since $\ch_qU^i\geq \ch_qL(c_i,0)$ for $i=1,2$ we
have $\ch_qL(c_1,0)\ch_q(c_2,0)\leq \ch_qV.$ So
$$\eta(q)\ch_qU^1\ch_qU^2\leq \eta(q)\ch_qV,\ \ 
\eta(q) \ch_qL(c_1,0)\ch_qL(c_2,0)\leq \eta(q)\ch_qV$$
 as functions for $q\in (0,1).$

Assume that $\ch_qU^1=\frac{q^{-c_1/24}}{\prod_{n>1}(1-q^n)}.$ 
Then 
$$\eta(q)\ch_qU\geq \eta(q)\frac{q^{-c_1/24}}{\prod_{n>1}(1-q^n)}
\ch_qL(c_2,0)$$
as functions for $q\in (0,1).$ That is,
$$\eta(q)\ch_qU\geq q^{c_2/24}(1-q)\ch_qL(c_2,0).$$ 
From the proof of Lemma \ref{l4.2} we see that if the coefficients of 
$(1-q)\ch_qL(c_2,0)$ satisfy the polynomial growth condition, so does
the coefficients of  $\ch_qL(c_2,0).$ But this is impossible by Lemma 
\ref{l4.2}. But $q^{c_2/24}(1-q)\ch_qL(c_2,0)\leq \eta(q)\ch_q V$ 
as functions for $q\in (0,1)$ and the coefficients of $\eta(q)\ch_q V$
satisfy the polynomial growth condition. This is a contradiction. 

By Proposition \ref{vir} we see immediately that $c_i=c_{s_i,t_i}$ for
some $s_i,t_i$ and $\<\omega_i\>$ is isomorphic to $L(c_{s_i,t_i},)$
for $i=1,2.$ \qed

\begin{lem}\label{l4.5} Let $c_i=c_{s_i,t_i}$ as in Lemma \ref{l4.4}. Then 
both $c_1$ and $c_2$  are $1/2.$
\end{lem}

\pf We need to solve the equation
$$1-\frac{6(s_1-t_1)}{s_1t_1}+1-\frac{6(s_2-t_2)}{s_2t_2}=1$$
for two pairs of coprime integers $1<s_i<t_i.$ That is,
$$\frac{s_1}{t_1}+\frac{t_1}{s_1}+\frac{s_2}{t_2}+\frac{t_2}{s_2}=\frac{25}{6}.$$
Let $x=\frac{s_1}{t_1}$ and $y=\frac{s_2}{t_2}.$ Then we the equation becomes
$$x+\frac{1}{x}+y+\frac{1}{y}=\frac{25}{6}.$$

The following argument using the elliptic curce is due to  N. Elkies and 
we thank him and A. Ryba for communicating the solution to us.
The equation $x+\frac{1}{x}+y+\frac{1}{y}=\frac{25}{6}$ gives 
an elliptic curve. Multiply the equation by $6xy$ to get
$$E: 6xy^2+6x^2y+6x+6y=25.$$ 
Putting one of the Weierstrass points at infinity yields the curve
$$y^2+XY=X^3-1070X+7812$$
which has rank 0 over $\Q.$ So every rational points
in $E$ is a torsion point. So  $E/\Q$ has at most 16 torsion points.
Note that the curve  has $8$ obvious symmetries, generated by the involutions
taking $(x,y)$ to $(1/x,y), (x,1/y),$ and $(y,x).$ Here are the 
rational points in $E:$ four from $(\frac{3}{4},\frac{3}{4}),$
four from $(1,\frac{2}{3}),$ four from  $(-1,6)$ and $4$ from infinity.

Since we assume that $1<s_i<t_i$ and $s_i,t_i$ are coprime, we immediately
see that the only solution intersting to us is $(\frac{3}{4},\frac{3}{4}).$
This is, $c_i=\frac{1}{2}$ for $i=1,2.$ \qed
\bigskip

Here is a characterization of $L(1/2,0)\otimes L(1/2,0).$
\begin{th}\label{t4.1} If $V$ is a rational and $C_2$-cofinite vertex operator
algebra of the moonshine type such that $c=\c=1$ and $\dim V_2=2,$ then $V$ is isomorphic to $L(1/2,0)\otimes L(1/2,0).$
\end{th}

\pf  By Lemmas \ref{l4.4} and \ref{l4.5}, the vertex operator subalgebra $U$ generated by $V_2$ of $V$ is isomorphic
to $L(\frac{1}{2},0)\otimes L(\frac{1}{2},0)$ which is rational
has 9 inequivalent irreducible modules $L(\frac{1}{2},h_1)\otimes L(\frac{1}{2},h_2)$ for $h_i\in \{0,\frac{1}{2},\frac{1}{16}\}$ (see 
\cite{DMZ},\cite{W}).Thus $V$ is a direct sum of irreducible 
$L(\frac{1}{2},0)\otimes L(\frac{1}{2},0)$-modules. Note that
$h_1+h_2\in \Z$ if and only if $h_1=h_2=0$ or
$h_1=h_2=\frac{1}{2}.$ So only   $L(\frac{1}{2},0)\otimes L(\frac{1}{2},0)$
and  $L(\frac{1}{2},\frac{1}{2})\otimes L(\frac{1}{2},\frac{1}{2})$
can possibly occur in $V$ as $L(\frac{1}{2},0)\otimes L(\frac{1}{2},0)$-modules.Using  the assumption that $\dim V_0=1$ and $V_1=0$ gives the result
that $V$ is isomorphic to $L(\frac{1}{2},0)\otimes L(\frac{1}{2},0).$
\qed
\bigskip

We end this paper with the following conjecture which 
strengthens Theorem \ref{t4.1}.
\begin{con}\label{cj} If $V$ is simple, rational and $C_2$ confinite vertex
operator algebra of the moonshine type with $c=\c=1$ and $\dim V_2>1,$
then $V$ is isomorphic to $L(1/2,0)\otimes L(1/2,0).$ 
\end{con}

It is essentially proved in \cite{K} that if $V$ is a rational vertex operator 
algebra such that $\sum_{i}|\chi_i(q)|^2$ is modular invariant where 
$\chi_i(q)$ are the $q$-character of the irreducible $V$-module, then
the $q$-character of $V$ is equal to the character of one of the following
vertex operator algebras $V_L, V_L^+$ and $V_{\Z\alpha}^G$ where
$L$ is any positive definite even lattice of rank 1, $V_L^+$ is the
fixed points of the automorphism of $V$ lifted from the
$-1$ isometry of $L,$ and $\Z\alpha$ is the root lattice
of type $A_1$ such that $(\alpha,\alpha)=2$ and $G$ is a finite subgroup
of $SO(3)$ isomorphic to $A_4,S_4$ or $A_5.$ It is widely believed that
 $V_L, V_L^+$ and $V_{\Z\alpha}^G$ should give a complete list of 
simple and rational vertex operator algebras with $c=\c=1.$ One can easily
find counter example if $c\ne \c.$ It is clearly from
the construction that if $V$ is one of these vertex operator algebras of the moonshine type then $\dim V_2= 2.$ This should be a very strong evidence for
the conjecture \ref{cj}. We remark that the assumption that  $\sum_{i}|\chi_i(q)|^2$ is modular invariant in \cite{K} is still an open problem in mathematics.

\end{document}